\documentclass[12pt]{amsproc}
\usepackage[ascii]{inputenc}
\usepackage[english]{babel}
\usepackage{amssymb}
\newtheorem*{tm}{Theorem}
\newtheorem*{tmP}{Pontrjagin Theorem}
\newtheorem*{tmL}{Langer Theorem}
\newtheorem{lem}{Lemma}
\newcommand{\rank}{\operatorname{rank}}
\renewcommand{\Re}{\operatorname{Re}}
\renewcommand{\Im}{\operatorname{Im}}

\title[On invariant subspaces]{On invariant subspaces of dissipative
operators in a space with indefinite metric}
\author{A.~A.~Shkalikov}
\thanks{This work is supported by the Russian Foundation for Basic
Research under grant No.~04-01-00412 and by the Foundation of Support
of Leading Scientific Schools under grant No.~1927.2003.1.}
\address{Department of Mechanics and Mathematics\\
Moscow Lomonosov State University\\
Moscow, 119992, Russia}
\email{ashkalikov@yahoo.com}
\begin{document}
\begin{abstract}
The theorem on the existence of maximal nonnegative invariant
subspaces for a special class of dissipative operators in Hilbert
space with indefinite inner product is proved in the paper. It is
shown in addition that the spectra of the restrictions of these
operators on the corresponding invariant subspaces lie in the
closed upper half-plane. The obtained theorem is a generalization
of well-known results of L.~S.~Pontrjagin, H.~K.~Langer, M.~G.~Krein
and T.~Ja.~Azizov devoted to this subject.
\end{abstract}
\begin{flushleft}
UDK~517.9+517.43
\end{flushleft}
\maketitle

\section{Introduction}
Let \(\mathcal H\) be a separable Hilbert space with usual scalar
product \((x,y)\) and indefinite one \([x,y]=(Jx,y)\), where
\(J=P_+-P_-\), and \(P_+\), \(P_-\) are the orthoprojectors
such that \(P_+P_-=P_-P_+=0\), \(P_++P_-=I\)
and \(I\) is the identity operator. Obviously, \(J\) admits
such a representation if and only if \(J=J^*\) and \(J^2=I\).
The space \(\{\mathcal H,J\}\) is called \emph{the Pontrjagin space}
and is denoted by \(\Pi_{\varkappa}\), if either \(\rank P_+\) or
\(\rank P_-\) is finite and equals \(\varkappa\). It is called
\emph{the Krein space} if both later numbers are infinite. A subspace
\(\mathcal L\) in \(\{\mathcal H, J\}\) is called \emph{nonnegative}
(\emph{uniformly positive}), if \([x,x]\geqslant 0\) (\(\geqslant
\varepsilon (x,x)\) with some \(\varepsilon\) independent on
\(x\)) for all \(x\in\mathcal L\). A nonnegative (uniformly positive)
subspace \(\mathcal L\) is said to be \emph{maximal} if there are
no nontrivial nonnegative (uniformly positive) extensions of this
subspace. Maximal nonpositive and uniformly negative subspaces are
defined analogiously.

Let us represent the space \(\mathcal H\) in the form
\(\mathcal H=\mathcal H^+\oplus\mathcal H^-\) where
\(\mathcal H^{\pm}=P_{\pm}(\mathcal H)\) are the ranges of the orthogonal
projectors \(P_{\pm}\). Consider a linear operator \(A\) in
\(\mathcal H\) with domain of definition \(\mathcal D(A)\). The spectrum
and the resolvent set of \(A\) is denoted further by \(\sigma(A)\)
and \(\rho(A)\). An operator \(A\) is called \emph{dissipative}
in \(\mathcal H\) if \(\Im(Ax,x)\geqslant 0\) for all \(x\in\mathcal D(A)\).
A dissipative operator is called \emph{maximal dissipative} if there are
no nontrivial dissipative extensions of this operator. It is
known~\cite[Ch.~V, \S\,3.10]{1} that the later condition holds if and
only if \(\rho(A)\supset\mathbb C^+\) where \(\mathbb C^+\) is the
open upper half-plane. An operator \(A\) is called
\emph{dissipative} (\emph{maximal dissipative}) in the space
\(\{\mathcal H,J\}\) if \(JA\) is dissipative (maximal dissipative)
in \(\mathcal H\). Analogiously, \(A\) is called \emph{symmetric}
(\emph{self-adjoint}) in the space \(\{\mathcal H, J\}\) if \(JA\)
is symmetric (self-adjoint) in the space \(\mathcal H\).

In the sequel, we work only with operators \(A\) for which the sum
\(\mathcal D^+\oplus\mathcal D^-\) is dense in \(\mathcal H\), where
\(\mathcal D^{\pm}=\mathcal D(A)\cap\mathcal H^{\pm}\). We will
always assume that \(\mathcal D(A)=\mathcal D^+\oplus\mathcal D^-\),
otherwise we can consider the restriction of \(A\) to this
domain. In this case the operator \(A\) can be represented as
an operator matrix with respect to the decomposition
\(\mathcal H=\mathcal H^+\oplus\mathcal H^-\):
\begin{equation}\label{eq:1}
	A=\begin{pmatrix}P_+AP_+&P_+AP_-\\ P_-AP_+&
	P_-AP_-\end{pmatrix}:=\begin{pmatrix}A_{11}&A_{12}\\
	A_{21}&A_{22}\end{pmatrix}.
\end{equation}
The vectors \(x=x_++x_-\in\mathcal H\) with \(x_{\pm}\in
\mathcal H_{\pm}\) are indentified in this representation with the
colomns \(x=\begin{pmatrix}x_+\\ x_-\end{pmatrix}\), and the action
of \(A\) is determined by the formula
\[
	Ax=A\begin{pmatrix}x_+\\ x_-\end{pmatrix}=
	\begin{pmatrix} A_{11}x_++A_{12}x_-\\
	A_{21}x_++A_{22}x_-\end{pmatrix},\qquad
	x_+\in\mathcal D^+,\; x_-\in\mathcal D^-.
\]

Pontrjagin~\cite{2} proved in~1944 the following fundamental result.

\begin{tmP}
Let \(A\) be a self-adjoint operator in the space \(\{\mathcal H,J\}\)
and \(\rank P_+=\varkappa<\infty\). Then there exists a maximal
nonnegative \(A\)-invariant subspace \(\mathcal L\)
(\(\dim\mathcal L=\varkappa\)) such that the spectrum of the
restriction \(A|_{\mathcal L}\) lies in the closed upper half-plane.
\end{tmP}

Starting from paper~\cite{2} the problem on the existence of maximal
definite invariant subspaces has been a key-stone of the operator theory
in Pontrjagin and Krein spaces. Krein~\cite{2a} obtained an analogue of
Pontrjagin theorem for unitary operators in \(\Pi_{\varkappa}\)
and developed a new approach to the problem in question. An important
generalization of Pontrjagin theorem was obtained by Langer~\cite{3,4}
and Krein~\cite{5}. Let us present here the result~\cite{4}.

\begin{tmL}
Let \(A\) be a selfadjoint operator in Krein space
\(\{\mathcal H,J\}\) and \(\mathcal D(A)\supset\mathcal H^+\)
(the later condition holds if and only if \(A\) admits
representation~\eqref{eq:1} where \(A_{11}\) and \(A_{12}\) are
bounded). If in addition the operator \(A_{12}=P_+AP_-\) is compact, then
there exists a maximal \(A\)-invariant subspace \(\mathcal L\) such that
the spectrum of the restriction \(A|_{\mathcal L}\) lies in the closed
upper half-plane.
\end{tmL}

Later on the theorems on the existence of \(A\)-invariant
subspaces have been obtained for other classes of operators.
Krein brought into consideration and investigated the class of
definite operators, and Langer~\cite{6,7} proved the theorem
on the existence of maximal definite invariant subspaces for a wider
class of the so-called definitizable operators and obtained for these
operators an analogue of the spectral theorem. Krein and
Langer~\cite{8} and independently Azizov~\cite{9}
showed that Pontrjagin theorem remains to be valid (as before in
Pontrjagin space \(\{\mathcal H,J\}\), \(\rank P_+=\varkappa<\infty\)))
if the condition for \(A\) to be self-adjoint is replaced by the
condition to be maximal dissipative. Later on, Azizov and
Khoroshavin~\cite{10} proved an analogue of Langer theorem for a class
of nonstretching operators in Krein space, and
Azizov~\cite[Ch.~2]{11} proved that Langer theorem~\cite{4} remains to
be valid for maximal dissipative operators in Krein space. A direct
and shorter proof of the later result was suggested by the author~\cite{12}.

The Langer condition \(\mathcal D(A)\supset\mathcal H^+\) (or equivalently
the boundedness of the operators \(A_{11}\), \(A_{21}\)) is rather
restrictive. In particular, often in concrete problems
(see~\cite{14,15}, for example) the operator \(A_{21}\) is unbounded.

\section{Main result}
The goal of the present paper is to obtain a generalization of
Pontrjagin--Krein--Langer--Azizov theorem dropping out the Langer condition
\(\mathcal D(A)\supset\mathcal H^+\), i.~e. the condition
for the operators \(A_{11}\) and \(A_{21}\) to be bounded. The essence
of the assumptions formulated below can be expressed as follows:
the operator \(A_{22}\) is dominant with respect to the
interlacing operators \(A_{21}\) and \(A_{12}\), and the so-called transfer-function
of the operator matrix~\eqref{eq:1} is bounded. Let us formulate the
main result.

\begin{tm}
Let \(A\) be a dissipative operator in Krein space
\(\{\mathcal H,J\}\) and its domain \(\mathcal D(A)=\mathcal D^+
\oplus\mathcal D^-\) be dense in \(\mathcal H=\mathcal H^+
\oplus\mathcal H^-\). Let~\eqref{eq:1} be the matrix representation of
\(A\) in \(\mathcal H^+\oplus\mathcal H^-\) and the following conditions
hold:
\begin{itemize}
\item[(i)] the operator \(-A_{22}\) is maximal dissipative in the space
\(\mathcal H^-\) (and hence the resolvent \((A_{22}-\mu)^{-1}\)
exists for all \(\mu\in\mathbb C^+\));
\item[(ii)] the operator \(F(\mu)=(A_{22}-\mu)^{-1}A_{21}\) admits
a bounded closure for some (and hence for all) \(\mu\in\mathbb C^+\);
\item[(iii)] the operator \(G(\mu)=A_{12}(A_{22}-\mu)^{-1}\) is compact
for some (and hence for all) \(\mu\in\mathbb C^+\);
\item[(iv)] the operator
\[
	S(\mu)=A_{11}-A_{12}(A_{22}-\mu)^{-1}A_{21}
\]
admits a bounded closure for some (and hence for all) \(\mu\in\mathbb C^+\).
\end{itemize}
Then the closure \(\overline{A}\) of the operator \(A\) is maximal dissipative
in the space \(\{\mathcal H,J\}\), and there exists a maximal nonnegative
\(\overline{A}\)-invariant subspace \(\mathcal L\) such that the spectrum of the
restriction \(\overline{A}|_{\mathcal L}\) lies in the closed upper half-plane.
Moreover, \(\mathcal L\subset\mathcal D(\overline{A})\), i.~e. the operator
\(\overline{A}|_{\mathcal L}\) is bounded.
\end{tm}

First we shall make two remarks on the conditions of the above theorem. It
is useful to view in mind that condition~(ii) is valid if the operator
\(A_{21}\) is closable (hence the adjoint operator \(A_{21}^*\) is densely
defined) and \(\mathcal D(A_{21}^*)\supset\mathcal D(A_{22}^*)\) (it is
known~\cite[Ch.~5]{1} that the adjoint to the dissipative operator
\(-A_{22}\) is densely defined). In fact, if the condition
\(\mathcal D(A_{21}^*)\supset\mathcal D(A_{22}^*)\) holds, then
\(F^*(\mu)=A_{21}^*(A_{22}^*-\overline{\mu})^{-1}\) is defined on the
whole \(\mathcal H^-\) and the adjoint to this operator is the closure
of the densely defined operator \(F(\mu)=(A_{22}-\mu)^{-1}A_{21}\).
Consequently, both operators \(F^*(\mu)\) and  \(\overline{F(\mu)}\)
are bounded. The second remark concerns condition~(i) which has not
been met in the formulations of the previous theorems on this subject.
However, it follows from~\cite[Ch.~2, Th.~2.9]{11} that if
\(\mathcal D(A)\supset\mathcal H^+\) then \(A\) is maximal dissipative
in \(\{\mathcal H, J\}\) if and only if \(-A_{22}\) is maximal dissipative
in \(\mathcal H^-\). Hence conditions~(i)--(iv) are weaker then those
in theorems of Pontrjagin, Krein, Langer and Azizov.

Later on, if we meet no confusions, we shall write \(\mu\) instead
of \(\mu I\) where \(I\) is the identity operator in \(\mathcal H^+\),
\(\mathcal H^-\) or in \(\mathcal H\).

\section{Preliminary propositions}
We shall premise several lemmas to the proof of Theorem.
Lemmas~\ref{lem:4} and~\ref{lem:5} play the key role.

\begin{lem}\label{lem:1}
A subspace \(\mathcal L\) is maximal nonnegative (uniformly positive)
if and only if it can be represented in the form
\begin{equation}\label{eq:2}
	\mathcal L=\left\{x=\begin{pmatrix}x_+\\ Kx_+\end{pmatrix},
	\qquad x_+\in\mathcal H^+\right\},
\end{equation}
where \(K:\mathcal H^+\to\mathcal H^-\) is a linear operator with the
norm \(\|K\|\leqslant 1\) (\(\|K\|<1\)). A nonnegative subspace
\(\mathcal L\) is maximal if and only if there exists no nonzero
element \(y_+\in\mathcal H^+\) such that \([x,y^+]=(x,y^+)=0\)
for all \(x\in\mathcal L\).
\end{lem}
\begin{proof}~(See~\cite{2}).
Assuming that \(\mathcal L\) is nonnegative subspace in
\(\{\mathcal H,J\}\) we have \(\|x_+\|\geqslant \|x_-\|\)
for all \(x=\begin{pmatrix}x_+\\x_-\end{pmatrix}
\in\mathcal L\). Then the restriction \(Q=P_+|_{\mathcal L}:
\mathcal L\to P_+(\mathcal L)\) is a bijection, and
\(\|Q^{-1}\|\leqslant 2\). Therefore,
\[
	\mathcal L=\left\{x=\begin{pmatrix}x_+\\ Kx_+\end{pmatrix},
	\qquad x_+\in P_+(\mathcal L),\qquad K=P_-Q^{-1}\right\}.
\]
Here \(\|K\|\leqslant 1\) if \(\mathcal L\) is nonnegative and
\(\|K\|< 1\) if \(\mathcal L\) is uniformly positive. Obviously,
\(\mathcal L\) is maximal if and only if \(P_+(\mathcal L)=
\mathcal H^+\). The second assertion of Lemma is also obvious.
\end{proof}

The operator \(K\) participating in representation~\eqref{eq:2}
is said to be \emph{the angle operator} of the subspace \(\mathcal L\).

\begin{lem}\label{lem:2}
Let \(A\) be an operator with dense domain \(\mathcal D(A)=
\mathcal D^+\oplus\mathcal D^-\), the resolvent set
\(\rho(A_{22})\) be nonempty, and the operators
\begin{equation}\label{eq:3}
	G=A_{12}(A_{22}-\mu)^{-1},\qquad F=(A_{22}-\mu)^{-1}A_{21},
	\qquad S=A_{11}-A_{12}F
\end{equation}
be bounded for some \(\mu\in\rho(A_{22})\). Then \(A\) is closable
and its closure is given by the relation
\begin{equation}\label{eq:4}
	\overline{A}=\mu+\begin{pmatrix}1&G\\0&1\end{pmatrix}
	\begin{pmatrix}S-\mu&0\\0&A_{22}-\mu\end{pmatrix}
	\begin{pmatrix}1&0\\F&1\end{pmatrix}.
\end{equation}
More precisely, the domain and the action of \(\overline{A}\)
are defined by the relations
\begin{gather*}
	\mathcal D(\overline{A})=\left\{\begin{pmatrix}x_+\\
	x_-\end{pmatrix}\in\mathcal H,\qquad x_+\in\mathcal H^+,
	\qquad Fx_++x_-\in\mathcal D^-\subset\mathcal D(A_{22})
	\right\},\\
	\overline{A}\begin{pmatrix}x_+\\x_-\end{pmatrix}=
	\begin{pmatrix}Sx_++G(A_{22}-\mu)(Fx_++x_-)\\
	(A_{22}-\mu)(Fx_++x_-)+\mu x_-\end{pmatrix}.
\end{gather*}
\end{lem}
\begin{proof}~(Cf.~\cite{14}).
One can easily check the validity of representation~\eqref{eq:4} for
\(x=\begin{pmatrix}x_+\\x_-\end{pmatrix}\in\mathcal D(A)\).
Since the operators \(G\), \(S\), \(F\) are bounded, we conclude that the first
and the third matrix in the right hand-side of~\eqref{eq:4}
are invertible, and the second one
represents a closed operator. Therefore, \(A\) is closable and
representation~\eqref{eq:4} is valid. The description of
\(\mathcal D(\overline{A})\) and the formula for the action of
\(\overline{A}\) follows from~\eqref{eq:4}.
\end{proof}

\begin{lem}\label{lem:3}
Suppose that \(-A_{22}\) is a maximal dissipative operator in
\(\mathcal H^-\) and \(G(\mu)=A_{12}(A_{22}-\mu)^{-1}\) is compact
for some \(\mu\in\mathbb C^+\). Then \(\|G(\mu)\|\to\infty\)
as \(\Im\mu\to +\infty\).
\end{lem}
\begin{proof}
It follows from the equation
\[
	G(\lambda)=G(\mu)+(\lambda-\mu)G(\mu)(A_{22}-\lambda)^{-1}
\]
that \(G(\lambda)\) is compact for all \(\lambda\in\mathbb C^+\).
Further we make use from the relation
\[
	G(\mu)=G(i)(A_{22}+i)(A_{22}-\mu)^{-1}.
\]
The norm of the operator function \((A_{22}+i)(A_{22}-\mu)^{-1}\)
is uniformly bounded in the half-plain \(\Im\mu\geqslant\varepsilon\).
The compact operator \(G(i):\mathcal H^-\to\mathcal H^+\) can be
approximated with arbitrary accuracy in the norm operator topology
by a finite rank operator. Hence it suffices to prove that
\(\|Q(A_{22}+i)(A_{22}-\mu)^{-1}\|\to 0\) as \(\Im\mu\to\infty\)
for any operator \(Q\) of rank \(1\), namely, for
\(Q=(\cdot,\varphi_-)\varphi_+\) where \(\varphi_{\pm}\in\mathcal H_{\pm}\).
Observe, that \(Q\) can be approximated with arbitrary accuracy in the
norm operator topology by an operator of the form \(Q_0=(\cdot,\varphi_0)
\varphi_+\) where \(\varphi_0\in\mathcal D(A_{22}^*)\) (we already noted that
the adjoint to a dissipative operator is densely defined). Now, the operator
\(Q_0(A_{22}+i)\) is bounded, and \(\|(A_{22}-\mu)^{-1}\|\leqslant
1/\Im\mu\) for \(\mu\in\mathbb C^+\). This gives the assertion of Lemma.
\end{proof}

\begin{lem}\label{lem:4}
Let the conditions of Lemma~\ref{lem:2} be preserved for an operator
\(\overline{A}\) as well as the notations~\eqref{eq:3} for the operators
\(G\), \(F\), \(S\), and \(\overline{A}\) for the closure of \(A\).
Then a subspace
\[
	\mathcal L=\left\{x=\begin{pmatrix}x_+\\ Kx_+\end{pmatrix},
	\qquad x_+\in\mathcal H^+\right\}
\]
with the angle operator \(K:\mathcal H^+\to\mathcal H^-\) lies in
\(\mathcal D(\overline{A})\) if and only if the operator
\(L=A_{21}+(A_{22}-\mu)K:\mathcal H^+\to\mathcal H^-\) is well
defined on \(\mathcal D^+\) and admits a bounded closure. If the later
condition holds, then the subspace \(\mathcal L\) is
\(\overline{A}\)-invariant if and only if
\begin{equation}\label{eq:5}
	L=K(S-\mu+GL),
\end{equation}
and then the restriction \(\overline{A}|_{\mathcal L}\) is
represented in the form
\begin{equation}\label{eq:5a}
	\overline{A}|_{\mathcal L}=Q^{-1}(S+GL)Q,
\end{equation}
where \(Q=P^+|_{\mathcal L}:\mathcal L\to\mathcal H^+\), \(\|Q^{-1}\|
\leqslant 1+\|K\|\).
\end{lem}
\begin{proof}
Let \(L=A_{21}+(A_{22}-\mu)K\) be defined on \(\mathcal D^+\) and admit a
bounded closure. Then
\[
	(A_{22}-\mu)^{-1}Lx_+=(F+K)x_+\in\mathcal D^-
\]
for all \(x_+\in\mathcal H^+\). Recalling the description of
\(\mathcal D(\overline{A})\) obtained in Lemma~\ref{lem:2} we find
\(\mathcal L\subset\mathcal D(\overline{A})\) and
\begin{equation}\label{eq:6}
	(\overline{A}-\mu)\begin{pmatrix}x_+\\ Kx_+\end{pmatrix}=
	\begin{pmatrix}(S-\mu+GL)x_+\\ Lx_+\end{pmatrix}.
\end{equation}
Conversely, if \(\mathcal L\subset\mathcal D(\overline{A})\), then
\((F+K)x_+\in\mathcal D^-\). Hence \(L=(A_{22}-\mu)(F+K)\) is defined
on the whole \(\mathcal H^+\). Since the operator \(\overline{A}:
\mathcal H\to\mathcal H\) is closed, its restriction \(\overline{A}:
\mathcal L\to\mathcal H\) is also closed. The later operator is defined
on the whole \(\mathcal L\). Therefore, it is bounded by virtue of the
closed graph theorem, and it follows from~\eqref{eq:6} that the operator
\(L\) is bounded. Now, suppose that the subspace \(\mathcal L\) is
\(\overline{A}\)-invariant. Then given \(x_+\in\mathcal H^+\) there
exists an element \(y_+\in\mathcal H^+\) such that
\begin{equation}\label{eq:7}
	\begin{pmatrix}(S-\mu+GL)x_+\\ Lx_+\end{pmatrix}=
	\begin{pmatrix}y_+\\ Ky_+\end{pmatrix}.
\end{equation}
This implies equation~\eqref{eq:5}. Conversely, suppose that \(\mathcal L
\subset\mathcal D(\overline{A})\) and equation~\eqref{eq:5} holds. Then
relation~\eqref{eq:7} is valid, and it is equivalent to
\(\overline{A}\)-invariance of the subspace \(\mathcal L\). The last assertion
of Lemma follows from~\eqref{eq:6}. We remark only that the estimates
\(\|Q\|\leqslant 1\) and \(\|Q^{-1}\|\leqslant 1+\|K\|\) follow from
the definition of \(Q\).
\end{proof}

\begin{lem}\label{lem:5}
Let \(A\) be an uniformly dissipative operator in the space
\(\{\mathcal H,J\}\), i.~e.
\begin{equation}\label{eq:8}
	\Im[Ax,x]\geqslant\varepsilon (x,x)\qquad\text{for }
	x\in\mathcal D(A),
\end{equation}
where \(\varepsilon>0\). Let \(\mathcal D(A)\supset\mathcal H^+\),
and \(-A_{22}\) be a maximal dissipative operator in
\(\mathcal H^-\). Then the operator \(A\) is maximal dissipative
in \(\{\mathcal H,J\}\), the real axis belongs to the resolvent set
\(\rho(A)\) and its spectrum in \(\mathbb C^+\) is bounded. If a
Jordan contour \(\Gamma_+\) surrounds the set \(\sigma(A)\cap
\mathbb C^+\) and
\begin{equation}\label{eq:9}
	Q_+=\dfrac{1}{2\pi i}\int\limits_{\Gamma_+}
	(\lambda-A)^{-1}\,d\lambda
\end{equation}
is the corresponding Riesz projector, then \(\mathcal L=Q_+(
\mathcal H)\) is an \(A\)-invariant maximal uniformly positive
subspace. Moreover,
\begin{equation}\label{eq:10}
	[x,x]\geqslant 2\varepsilon(\pi\|A_+\|)^{-1}\,\|x\|^2,
	\qquad\text{for }x\in\mathcal L
\end{equation}
where \(A_+=A|_{\mathcal L}\). If for some \(\mu\in\mathbb C^+\)
the estimate \(\|G(\mu)\|=\gamma<1\) holds, then
\begin{equation}\label{eq:11}
	\|A_+\|\leqslant 2(\|S\|+\gamma(1-\gamma)^{-1}
	(\|S\|+|\mu|)),\qquad S=A_{11}-G(\mu)A_{21}.
\end{equation}
\end{lem}
\begin{proof}
We already noted reffering to~\cite[Ch.~2, \S\,2]{11} that \(A\) is
maximal dissipative in \(\{\mathcal H,J\}\) under the assumptions of this
Lemma (in fact, if \(JA\) admits a nontrivial dissipative extension in
\(\mathcal H\), then the condition
\(\mathcal D(A)\supset\mathcal H^+\) implies that \(-A_{22}\)
admits nontrivial dissipative extensions in \(\mathcal H^-\), and we
come to a contradiction). The other assertions of Lemma but the last
estimates were proved in the author's paper~\cite{12}. To prove
estimate~\eqref{eq:11} we have to repeat partially the arguments
from~\cite{12}. We do this in several steps.

\textit{Step~1.} It follows from the condition
\(\mathcal D(A)\supset\mathcal H^+\) that the operator \(AP_+\)
is bounded. Take a number \(a>2\|AP_+\|\). Denote as before
\(G(\lambda)=A_{12}(A_{22}-\lambda)^{-1}\) and show that
\begin{equation}\label{eq:12}
	\|G(\lambda)\|\leqslant 2+2a\varepsilon^{-1}\qquad
	\text{for all }\lambda\in\overline{\mathbb C^+}.
\end{equation}
Consider the operator
\[
	T(\lambda,a)=\begin{pmatrix}ia&A_{12}\\
	0&-A_{22}+\lambda\end{pmatrix}=(JA+ia)+
	(\lambda-ia)P_--AP_+.
\]
The operators \(JA+ia\) and \(T(\lambda,a)+AP_+\) are
maximal dissipative for \(\lambda\in\mathbb C^+\). Moreover, we have the
estimate
\[
	\Im(T(\lambda,a)x,x)\geqslant(a/2)\|x\|^2\qquad
	\text{for }x\in\mathcal D(T)=\mathcal H^+\oplus
	\mathcal D(A_{22}),
\]
provided that \(\lambda\in\mathbb C^+_a=\{\lambda\mid\Im\lambda\geqslant a\}\).
Therefore \(T(\lambda,a)\) is invertible for \(\lambda\in\mathbb C^+_a\) and
\(\|T^{-1}(\lambda,a)\|\leqslant 2a^{-1}\). Since \(G(\lambda)=
aP_+T^{-1}(\lambda,a)P_-\), we have \(\|G(\lambda)\|\leqslant 2\) for
\(\lambda\in\mathbb C^+_a\). Using the equation
\[
	G(\lambda)=G(\lambda+ia)+iaG(\lambda+ia)(A_{22}-\lambda)^{-1}
\]
we get estimate~\eqref{eq:12}. Here we view in mind that
\(\Im(A_{22}x,x)\leqslant-\varepsilon (x,x)\) and hence
\(\|(A_{22}-\lambda)^{-1}\|\leqslant\varepsilon^{-1}\) for
\(\lambda\in\overline{\mathbb C^+}\).

\textit{Step~2.} It follows from representation~\eqref{eq:4} that
\(\lambda\in\rho(A)\cap\mathbb C^+\) if and only if the operator
\(S(\lambda)-\lambda\) is invertible. Since the operator \(A_{11}\)
is bounded, we have \((A_{11}-\lambda)^{-1}=-\lambda^{-1}+O(\lambda^{-2})\)
as \(\lambda\to\infty\). Viewing in mind that \(A_{21}\) is bounded and
 \(G(\lambda)\) is subject to estimate~\eqref{eq:12} in \(\mathbb C^+\)
we find
\begin{equation}\label{eq:13}
	(S(\lambda)-\lambda)^{-1}=(A_{11}-\lambda)^{-1}(1-G(\lambda)
	A_{21}(A_{11}-\lambda)^{-1})^{-1}=-\lambda^{-1}+O(\lambda^{-2}),
\end{equation}
as \(\lambda\in\overline{\mathbb C^+}\) and \(\lambda\to\infty\). Hence the
spectrum of \(A\) in \(\mathbb C^+\) is bounded.

\textit{Step~3.} Take a contour \(\Gamma_+\) consisting of a segment
\([-R,R]\) and the half of the circle \(C_R\) in \(\mathbb C^+\) of the
radius \(R\) and the center at zero. Taking \(R\) sufficiently large we
may insure that the spectrum of \(A\) in \(\mathbb C^+\) lies inside
\(\Gamma_+\). Consider the Riesz projector~\eqref{eq:9}. Obviously, the
subspace \(\mathcal L=Q_+(\mathcal H)\) is \(A\)-invariant, and the
restriction \(A_+=A|_{\mathcal L}\) is a bounded operator. We can replace
\(A\) by \(A_+\) in~\eqref{eq:9}. Then we have as \(R\to\infty\)
\[
	\dfrac{1}{2\pi i}\int\limits_{C_R}(\lambda-A_+)^{-1}\,
	d\lambda=\dfrac{1}{2\pi i}\int\limits_{C_R}(\lambda^{-1}+
	O(\lambda^{-2}))^{-1}\,d\lambda=\dfrac12 I+O(R^{-1}).
\]
Let \(x=Q_+x\in\mathcal L\) and \(y=(\lambda-A_+)^{-1}x\). Then
\[
	[x,x]=\Re[Q_+x,x]=\dfrac12[x,x]+\dfrac{1}{2\pi}
	\int\limits_{-R}^R\Im[y,(\lambda-A)y]\,d\lambda+
	O(R^{-1}).
\]
For \(\lambda\in\mathbb R\)
\[
	\Im[y,(\lambda-A)y]=\Im[Ay,y]\geqslant\varepsilon(y,y).
\]
Passing to the limit as \(R\to\infty\) and taking into account
the inequality
\[
	\|x\|\leqslant\|\lambda-A_+\|\,\|y\|\leqslant
	(|\lambda|+\|A_+\|)\,\|y\|,
\]
we get
\[
	[x,x]\geqslant\dfrac{1}{\pi}\int\limits_{-\infty}^{\infty}
	\Im[Ay,y]\,d\lambda\geqslant\dfrac{\varepsilon}{\pi}
	\|x\|^2\,\int\limits_{-\infty}^{\infty}(\|A_+\|+
	|\lambda|)^{-2}\,d\lambda=\dfrac{2\varepsilon}{\pi\|A_+\|}
	\|x\|^2.
\]
This proves that \(\mathcal L\) is uniformly positive subspace, and the
estimate~\eqref{eq:10} is valid.

\textit{Step~4.} Let us prove that \(\mathcal L\) is a maximal uniformly
positive subspace. It easily follows from~\eqref{eq:4} that
\[
	(\lambda-A)^{-1}=\begin{pmatrix}(\lambda-S(\lambda))^{-1}&*\\
	*&*\end{pmatrix},
\]
where by \(*\) we assign operators which representation is not used in the
sequel. For \(z\in\mathbb H^+\) we have
\begin{equation}\label{eq:14}
	((\lambda-A)^{-1}z,z)=((\lambda-S(\lambda))^{-1}z,z).
\end{equation}
Integrating the function \((2\pi i)^{-1}((\lambda-A)^{-1}z,z)\) along the
contour \(\Gamma_+\), using relations~\eqref{eq:13} and~\eqref{eq:14} for
the integrals along the half circle \(C_R\) and passing to the limit as
\(R\to\infty\), we obtain
\[
	(Q_+z,z)=[Q_+z,z]=\dfrac12(z,z)+\dfrac{1}{2\pi}
	\int\limits_{-\infty}^{\infty}\Im[Ay,y]\,d\lambda,
\]
where \(y=(\lambda-A)^{-1}z\). Consequently, \(2(Q_+z,z)\geqslant(z,z)\)
for all \(z\in\mathcal H^+\). Hence there is no nonzero element
\(z\in\mathcal H^+\) such that \(z\perp\mathcal L=Q_+(\mathcal H)\).
By lemma~\ref{lem:1} the subspace \(\mathcal L\) is maximal positive.

\textit{Step~5.} Finally, let us prove estimate~\eqref{eq:11} provided
that \(\|G(\mu)\|=\gamma<1\) for some \(\mu\in\mathbb C^+\). Since
\(\mathcal L\subset\mathcal D(A)\) and \(\mathcal L\) is an
\(A\)-invariant subspace, by Lemma~\ref{lem:3} we have
\[
	L=K(S-\mu+GL),
\]
where \(K\) is the angle operator of the subspace \(\mathcal L\),
and \(L=A_{21}+(A_{22}-\mu)K\). Consequently,
\[
	L=(1-KG)^{-1}K(S-\mu),\qquad \|L\|\leqslant(1-\gamma)^{-1}
	(\|S\|+|\mu|).
\]
From~\eqref{eq:5a} we get the inequality
\[
	\|A_+\|\leqslant 2(\|S\|+\gamma\|L\|),
\]
which implies estimate~\eqref{eq:11}. Lemma is proved.
\end{proof}

\begin{lem}\label{lem:6}
Let a sequence of linear operators \(T_n\) in the space \(\mathcal H\)
converge in the norm operator topology to an operator \(T\). Suppose that
the spectrum of \(T\) in a domain \(\Omega\subset\mathbb C\) is discrete.
If \(\sigma(T_n)\cap\Omega=\varnothing\) for all \(n\), then
\(\sigma(T)\cap\Omega=\varnothing\).
\end{lem}
\begin{proof}
Given \(\mu\in\Omega\cap\rho(T)\) there exists a neighbourhood
\(\mathcal U_{\delta}(\mu)\) such that \(\mathcal U_{\delta}(\mu)
\subset\rho(T)\cap\rho(T_n)\) for all sufficiently large \(n\),
and
\begin{equation}\label{eq:15}
	\|(T_n-\lambda)^{-1}-(T-\lambda)^{-1}\|\to 0\qquad
	\text{as }n\to\infty
\end{equation}
uniformly for \(\lambda\in\mathcal U_{\delta}(\mu)\). Take an arbitrary
contour \(\Gamma\) in \(\Omega\) which does not intersect the discrete
spectrum of \(T\). Taking from the cover \(\{\mathcal U_{\delta}(\mu)
\}_{\mu\in\Gamma}\) a finite subcover, using relation~\eqref{eq:15}
and viewing in mind that the spectra of the operators \(T_n\) are
empty inside \(\Gamma\), we obtain that the Riesz projector of
\(T\) along the contour \(\Gamma\) equals zero (since the
corresponding Riesz projectors of \(T_n\) equal zero). Consequently,
the spectrum of \(T\) inside \(\Gamma\) is empty. By arbitrary
choice of \(\Gamma\) the same is true inside \(\Omega\).
\end{proof}

\section{Proof of Theorem}
Take a system of linear independent elements \(\{\varphi_n\}_1^{\infty}\)
belonging to \(\mathcal D^+=\mathcal H^+\cap\mathcal D(A)\) such that
the linear span of this system is dense in \(\mathcal H^+\).
Denote by \(\mathcal H_n^+\) the finite dimensional subspaces with the
bases \(\{\varphi\}_1^n\) and by \(P_n\) the orthoprojectors on these
subspaces. Consider the operator
\[
	A_{n,\varepsilon}=\begin{pmatrix}P_nA_{11}P_n&P_nA_{12}\\
	A_{21}P_n&A_{22}\end{pmatrix}+i\varepsilon J,\qquad
	\varepsilon>0,
\]
acting in the space \(\mathcal H_n=\mathcal H_n^+\oplus\mathcal H^-\)
with the domain \(\mathcal D(A_{n,\varepsilon})=\mathcal H_n^+
\oplus\mathcal D^-\subset\mathcal D(A)=\mathcal D^+\oplus\mathcal D^-\).

Let us sketch the plan of the proof. The conditions of Lemma~\ref{lem:5}
are fulfiled for the operators \(A_{n,\varepsilon}\), since
\[
	\Im [A_{n,\varepsilon}x,x]=\Im [Ax,x]+\varepsilon(x,x)
	\qquad\text{for } x\in\mathcal D(A_{n,\varepsilon}).
\]
By virtue of Lemmas~\ref{lem:4} and~\ref{lem:5} there exist maximal
uniformly positive subspaces \(\mathcal L_n\) with the angle operators
\(K_n:\mathcal H_n^+\to\mathcal H^-\), \(\|K_n\|<1\), such that
\begin{equation}\label{eq:16}
	(1-K_nG)L_n=K_n(S_n-\mu),
\end{equation}
where \(\mu\in\mathbb C^+\) and
\[
	G=G(i\varepsilon+\mu),\qquad
	S_n=P_n(i\varepsilon+S(i\varepsilon+\mu))P_n,\qquad
	L_n=A_{21}P_n+(A_{22}-i\varepsilon-\mu)K_n.
\]
We remark that we can write \(G\) in equation~\eqref{eq:16}
instead of \(G_n=P_nG\), since \(K_nG_n=K_nG\).

It will be shown that one can pass to the limit in the weak operator
topology in equation~\eqref{eq:16} choosing a subsequence
\(n_k\to\infty\). The limit equation
\begin{equation}\label{eq:17a}
	(1-KG)L=K(S-\mu),
\end{equation}
holds with \(\|K\|<1\), \(L=A_{21}+(A_{22}-i\varepsilon-\mu)K\),
\(\|L\|\leqslant\mathrm{const}\). By virtue of Lemmas~\ref{lem:1}
and~\ref{lem:4} the subspace \(\mathcal L\) with the angle operator
\(K\) is \(\overline{A}+i\varepsilon J\)-invariant and maximal uniformly
positive. We remark that one can hardly realize a direct proof of the
analogue of Lemma~\ref{lem:5} for operators of the form
\(\overline{A}+i\varepsilon J\), \(\varepsilon>0\), since there is no
simple way to get representation~\eqref{eq:14} for \(S(\lambda)-\lambda\).
Further, the operators \(K\), \(G\), \(L\), and \(S\) in equation~\eqref{eq:17a}
depend on \(\varepsilon\). Choosing a proper subsequence \(\varepsilon_n\to 0\)
one can pass again to the limit in the weak operator topology and obtain
equation~\eqref{eq:17a} with an operator \(K\), \(\|K\|\leqslant 1\), and the
operators \(L=A_{21}+(A_{22}-\mu)K\) and \(S=S(\mu)\). Here \(L\) is bounded
and by Lemmas~\ref{lem:1} and \ref{lem:4} the subspace \(\mathcal L\)
with the angle operator \(K\) is \(\overline{A}\)-invariant and maximal
nonpositive. On this way we have also to prove that the spectra of the
restrictions \(\overline{A}+i\varepsilon J\) onto the invariant subspaces
\(\mathcal L_{\varepsilon}\) lie in the upper half-plane \(\mathbb C^+\)
for \(\varepsilon>0\) and in \(\overline{\mathbb C^+}\) for \(\varepsilon=0\).
From now on we realize the above plan.

By virtue of Lemma~\ref{lem:3} we can choose a number \(\mu\in\mathbb C^+\)
such that \(\|G\|=\|G(i\varepsilon+\mu)\|<1/2\) for all \(0\leqslant
\varepsilon\leqslant 1\). The operator function \(i\varepsilon+
S(i\varepsilon+\mu)\) is continuous for \(0\leqslant\varepsilon
\leqslant 1\) in the norm operator topology. Hence there is a
constant \(c\) such that
\begin{equation}\label{eq:17}
	\|i\varepsilon+S(i\varepsilon+\mu)\|\leqslant c\qquad
	\text{for all }0\leqslant\varepsilon\leqslant 1.
\end{equation}
It follows from~\eqref{eq:16} that \(\|L_n\|\leqslant 2(c+|\mu|)\).
We remark that
\[
	[x,x]\geqslant\delta(x,x)\qquad\text{for all }
	x\in\mathcal L_n
\]
if and only if \(\|K_n\|\leqslant 1-\delta\), where \(K_n\) is the angle
operator of the subspace \(\mathcal L_n\). By Lemma~\ref{lem:5} there is
a number \(\delta>0\) such that \(\|K_n\|\leqslant 1-\delta\). The operators
\(K_n\) and \(L_n\) acting from \(\mathcal H_n^+\) into \(\mathcal H^-\)
can be treated as operators from \(\mathcal H^+\) into \(\mathcal H^-\)
after their zero extension on the orthogonal complement \(\mathcal H^+
\ominus\mathcal H_n^+\). Certainly, the norms of these operators are
preserved. Since \(\mathcal H^+\) and \(\mathcal H^-\) are separable
spaces and \(\|K_n\|<1-\delta\), one can choose a weakly convergent subsequence
\(K_{n_j}\rightharpoonup K\) (here we make use of the fact that the unit
ball of a separable Hilbert space is a compact set in the weak topology).
Since the norms of the operators \(\{L_{n_j}\}\) are bounded by a
constant \(2(c+|\mu|)\), one can choose from the sequence \(\{L_{n_j}\}\)
a weakly convergent subsequence. Hence there are indices \(m=n_k\to\infty\)
such that \(K_m\rightharpoonup K\), \(L_m\rightharpoonup L\). Let us prove that
\begin{equation}\label{eq:18}
	L = A_{21}+(A_{22}-i\varepsilon-\mu)K.
\end{equation}
We have \((A_{22}-i\varepsilon-\mu)^{-1}L_m=F(\mu+i\varepsilon)P_m+K_m
\rightharpoonup F(\mu+i\varepsilon)+K\). Consequently, \((A_{22}-i\varepsilon
-\mu)^{-1}L=F(\mu+i\varepsilon)+K\), and this implies relation~\eqref{eq:18}.

Now we remark that the weak convergence \(K_n\rightharpoonup K\) implies
\(K_nG\rightharpoonup KG\) and \(GK_n\rightharpoonup GK\) for any bounded
operator \(G\). One can not guarantee the convergence
\(K_nGL_n\rightharpoonup KGL\) provided that the sequences
\(\{K_n\}\) and \(\{L_n\}\) are weakly convergent. However, the convergence
\(K_nGL_n\rightharpoonup KGL\) does hold if \(G\) is a compact operator (in
this case the convergence holds even in the norm operator topology). In fact,
a compact operator \(G\) can be approximated with arbitrary accuracy in the
norm operator topology by a finite rank operator, therefore it suffices to
prove the convergence for an operator \(G=(\cdot,v)u\) of rank \(1\). In
the later case we have for all \(x\in\mathcal H^+\) and \(y\in
\mathcal H^-\)
\[
	(K_nGL_nx,y)=(L_nx,v)(K_nu,y)\to(Lx,v)(Ku,y)=
	(KGLx,y).
\]
Hence, one can pass to the weak limit in~\eqref{eq:16} and obtain the
relation
\begin{equation}\label{eq:22}
	(1-KG)L=K(S(i\varepsilon+\mu)+i\varepsilon-\mu),
\end{equation}
where \(L=A_{21}+(A_{22}-i\varepsilon-\mu)K\) is a bounded operator and
\(\|K\|\leqslant 1-\delta\) with some \(\delta>0\). As we already mentioned,
by Lemmas~\ref{lem:1} and~\ref{lem:4} the subspace \(\mathcal L\) with the
angle operator \(K\) is \(\overline{A}+i\varepsilon J\)-invariant and maximal
uniformly positive. The restriction \(A_{\varepsilon}^+=
(A+i\varepsilon J)|_{\mathcal L}\) is a bounded uniformly dissipative operator
on the subspace \(\mathcal L\) with the inner product \([\;,\;]\), which is
equivalent to the usual inner product in \(\mathcal L\), since the subspace
\(\mathcal L\) is uniformly positive. Consequently the spectrum of this
restriction lies in the open upper half-plane \(\mathbb C^+\).

Now, we shall pass to the limit choosing a subsequence \(\varepsilon_n\to 0\).
Observe that
\begin{equation}\label{eq:21}
	\begin{gathered}
	G(\mu+i\varepsilon)=G(\mu)+i\varepsilon G(\mu)(A_{22}-
	i\varepsilon-\mu)^{-1}\\
	S(\mu+i\varepsilon)=S(\mu)+i\varepsilon G(\mu+i\varepsilon)F\mu).
	\end{gathered}
\end{equation}
Since \(\|G(\mu+i\varepsilon)\|<1/2\) for \(\varepsilon\geqslant 0\),
it follows from~\eqref{eq:22} that
\[
	\|L(\varepsilon)\|\leqslant 2(\|S(\mu)\|+\varepsilon
	(1+\|G(\mu)\|\,\|F(\mu)\|)+\mu),
\]
i.~e. the norms of \(L(\varepsilon)\) are uniformly bounded for
\(0<\varepsilon\leqslant 1\). Take any sequence \(K=K(\varepsilon_n)\)
and choose a weakly convergent subsequence \(K(\varepsilon_{n_j})\).
Further, choose a weakly convergent subsequence from the sequence
\(L(\varepsilon_{n_j})\). On this way we find numbers \(\varepsilon_m\to 0\)
such that \(K_n=K(\varepsilon_m)\rightharpoonup K\), \(L_n=L(\varepsilon_n)
\rightharpoonup L\). We can repeat the arguments applied while making the
first limit procedure and obtain the relation \(L(\mu)=A_{21}+(A_{22}-\mu)K\).

Taking into account relations~\eqref{eq:21} and recalling that the operator
\(G(\mu)\) is compact we can pass to the weak limit in relation~\eqref{eq:22}
as \(\varepsilon_m\to 0\). Thus we obtain that relation~\eqref{eq:22} holds
with \(\varepsilon =0\) and the operators \(K\), \(L\), \(\|K\|\leqslant 1\),
\(\|L\|\leqslant\mathrm{const}\). By Lemma~\ref{lem:4} the subspace
\(\mathcal L\) with the angle operator \(K\) is \(\overline{A}\)-invariant
and maximal nonnegative. From Lemma~\ref{lem:4} we also have
\[
	\overline{A}|_{\mathcal L}=Q^{-1}(S(\mu)+G(\mu)L(\mu))Q.
\]
It was already proved that the spectra of the operators
\[
	T(\varepsilon)=S(\mu+i\varepsilon)+G(\mu+i\varepsilon)
	L(\mu+i\varepsilon)+i\varepsilon
\]
lie in \(\mathbb C^+\) for each \(\varepsilon>0\). It follows from
relations~\eqref{eq:21} that \(T(\varepsilon)-i\varepsilon=T(0)+C\),
where \(C\) is a compact operator. Hence the spectrum of \(T(0)\) in the
half-plane \(\Im\lambda\geqslant-\varepsilon\) is discrete. Here
\(\varepsilon>0\) is arbitrary number, therefore spectrum of \(T(0)\)
in the open lower-half plane is discrete. From~\eqref{eq:21} we obtain
\(T(\varepsilon_n)\Rightarrow T(0)\) taking into account that
\(K_nGL_n\Rightarrow KGL\) if \(G\) is a compact operator. By
Lemma~\ref{lem:6} the spectrum of \(T(0)\) (and hence the spectrum of
\(\overline{A}|_{\mathcal L}\)) lies in \(\overline{\mathbb C^+}\).

It is left to prove that \(\overline{A}\) is a maximal dissipative
operator in the space \(\{\mathcal H,J\}\). It follows from~\eqref{eq:4}
that
\[
	\overline{A}-\mu+i\alpha P_+=\begin{pmatrix}1&G(\mu)\\
	0&1\end{pmatrix}\begin{pmatrix}S(\mu)-\mu+i\alpha&0\\
	0&0\end{pmatrix}\begin{pmatrix}1&0\\F(\mu)&1\end{pmatrix},
\]
provided that \(\mu\in\mathbb C^+\). Here the number \(\alpha>\Im\mu\)
can be choosen sufficiently large to guarantee the invertibility of the
operator \(S(\mu)-\mu+i\alpha\). In this case the operator
\(J(\overline{A}-\mu+i\alpha)\) is dissipative in \(\mathcal H\) and
invertible. Therefore the dissipative operator \(J\overline{A}\) is maximal
dissipative. This ends the proof of Theorem.

\end{document}